\documentclass[11pt]{article}
\usepackage{amsthm}
\usepackage{latexsym, amssymb, amsmath}

\begin{document}
\numberwithin{equation}{section}

\newtheorem{THEOREM}{Theorem}
\newtheorem{PRO}{Proposition}
\newtheorem{XXXX}{\underline{Theorem}}
\newtheorem{CLAIM}{Claim}
\newtheorem{COR}{Corollary}
\newtheorem{LEMMA}{Lemma}
\newtheorem{REM}{Remark}
\newtheorem{EX}{Example}
\newenvironment{PROOF}{{\bf Proof}.}{{\ \vrule height7pt width4pt depth1pt} \par \vspace{2ex} }
\newcommand{\Bibitem}[1]{\bibitem{#1} \ifnum\thelabelflag=1 
  \marginpar{\vspace{0.6\baselineskip}\hspace{-1.08\textwidth}\fbox{\rm#1}}
  \fi}
\newcounter{labelflag} \setcounter{labelflag}{0}
\newcommand{\labelon}{\setcounter{labelflag}{1}}
\newcommand{\Label}[1]{\label{#1} \ifnum\thelabelflag=1 
  \ifmmode  \makebox[0in][l]{\qquad\fbox{\rm#1}}
  \else\marginpar{\vspace{0.7\baselineskip}\hspace{-1.15\textwidth}\fbox{\rm#1}}
  \fi \fi}

\newcommand{\LEFTLINE}{\ifhmode\newline\else\noindent\fi}
\newcommand{\RIGHTLINE}[1]{\LEFTLINE\rightline{#1}}
\newcommand{\CENTERLINE}[1]{\LEFTLINE\centerline{#1}}
\def\BOX #1 #2 {\framebox[#1in]{\parbox{#1in}{\vspace{#2in}}}}
\parskip=8pt plus 2pt
\def\AUTHOR#1{\author{#1} \maketitle}
\def\Title#1{\begin{center}  \Large\bf #1 \end{center}  \vskip 1ex }
\def\Author#1{\vspace*{-2ex}\begin{center} #1 \end{center}  
 \vskip 2ex \par}
\renewcommand{\theequation}{\arabic{section}.\arabic{equation}}
\def\bdk#1{\makebox[0pt][l]{#1}\hspace*{0.03ex}\makebox[0pt][l]{#1}\hspace*{0.03ex}\makebox[0pt][l]{#1}\hspace*{0.03ex}\makebox[0pt][l]{#1}\mbox{#1} }
\def\psbx#1 #2 {\mbox{\psfig{file=#1,height=#2}}}

 
 
\Title{Nonnegative Trigonometric Polynomials,\\ Sturm's Theorem, and Symbolic Computation}

\begin{center}
MAN KAM KWONG\footnote{The research of this author is supported by the Hong Kong Government GRF Grant PolyU 5012/10P and the Hong Kong Polytechnic University Grants G-YK49 and G-U751}
\end{center}

\begin{center}
\emph{Department of Applied Mathematics\\ The Hong Kong Polytechnic University,\\ Hunghom, Hong Kong}\\
\tt{mankwong@polyu.edu.hk}\\[4ex]
February 27, 2014 \\
Revised: April 26, 2016
\end{center}

\par\vspace*{\baselineskip}\par

\newcommand{\mb}{\mathbf}
\newcommand{\Cr}{\color{red}}

\begin{abstract}
\parskip=6pt
In this paper, we explain a procedure based on a classical result of Sturm that 
can be used to determine rigorously whether a
given trigonometric polynomial is nonnegative in a certain interval or not.
Many examples are given.
This technique has been employed by the author in several recent works.

The procedure often involves tedious computations that are 
time-consuming and error-prone. Fortunately, symbolic
computation software is available to automate the procedure.
In this paper, we give the details of its implementation
in MAPLE 13. Some who are strongly attached to a more traditional
theoretical research framework may find such details boring or 
even consider
computer-assisted proofs suspicious. However, we emphasize again that the
procedure is completely mathematically rigorous.

\noindent
\underline{$\vphantom{y}$Addendum}. The author found out recently that the Sturm procedure has
been used previously to study nonnegative cosine polynomials by Brown, 
Wang and Wilson \cite{BWW}. Borwn, Koumandos and Wang \cite{bkw} 
have also used the same
technique to study sums of Legendre polynomials.
\end{abstract}

{\bf{Mathematics Subject Classification (2010).}} 26D05, 42A05.

{\bf{Keywords.}} Trigonometric sums, positivity, inequalities, symbolic
computation.

\newpage

\section{Introduction}
For convenience, we use the following acronyms:
AP (algebraic polynomial(s)), CP (cosine polynomial), NN (nonnegative), SP (sine polynomial), TP (trigonometric polynomial).

In this paper, a finite TP (trigonometric polynomial) of degree $ n $ is a sum of the form
\begin{equation}  a_0 + \sum_{k=1}^{n} (a_k\,\cos(kx) + b_k\,\sin(kx)) ,  \end{equation}
where $ a_k $ and $ b_k $ are real numbers. When all the $b_k=0$ ($a_k=0$),
the sum is a CP (cosine polynomial) (SP (sine polynomial)).

The study of TP that are NN (nonnegative) in a certain interval, typically $ [0,2\pi ] $ or $ [0,\pi ] $,
has a long history and has many applications.
Many authors have written excellent surveys on the topic. Please refer to
\cite{A}--\cite{Ks2}, \cite{rs} and the references therein.

If a TP is NN in $ [0,2\pi ] $, then it is NN for all real $ x $.
A necessary and sufficient condition in
terms of the positive definiteness of a certain Toeplitz matrix is known.
For methods to construct examples, see \cite{D}.

For the interval $ [0,\pi ] $ or other subintervals of $ [0,2\pi ] $, there are
no analogous necessary and sufficient conditions.
Many general sufficient criteria have been obtained. Among them are
the famous result of V.L.\ Vietoris and its extension due to A.S.\ Belov 
(see, for example, \cite{Ks2}).
For SP, under the additional assumption that $ b_k $ are monotonically
decreasing in $ k $, Belov's criterion is both necessary and sufficient, while the corresponding criterion for
CP is only sufficient. Without the monotonicity assumption, few
general results are known. Recent work by the author has made some progress
in this direction.

In this paper, we present a procedure that can be used to determine
rigorously whether a given TP, with specific coefficients,
is NN in a given interval or not.
The procedure is based on the classical Sturm Theorem that gives the
exact number of real roots of an algebraic polynomial in a given interval.
The procedure always works when the coefficients are all rational numbers. 
In case some or all of the coefficients are irrational, or if at one of the 
endpoints of the interval, $ \cos(x) $ is difficult to compute analytically,
the procedure may not be practical. A modified procedure using
numerical approximation is needed.

One may question the usefulness of such a procedure.
A good NN criterion should be applicable to a wide class, preferably 
including polynomials of all degrees. If the procedure is only able to settle
one specific polynomial at a time, the best we can get out of it is a finite
sample of particular NN polynomials. It will never lead us to any useful general result.

It has happened more than once in our quest for a general criterion that the arguments
we devised work well for polynomials of high degrees, but not for a 
finite collection of small values of $ n $. The Sturm procedure can then be 
used to complete the proof. In addition, the proof of a general result may require
a few trigonometric inequalities that can be proved using the Sturm procedure. See Section~7
for some concrete examples.

Although the procedure is theoretically rigorous, the amount of computation
involved can be tedious, time-consuming and error-prone. In Section~5, we
describe how modern symbolic computation software (in our investigation, we
have used MAPLE 13) can come to our aid. Indeed, the software is so
convenient to use that the author has exploited it to experiment with large
numbers of TP to look for patterns and formulate possible conjectures 
which later are confirmed.
Details of the software
usage are given so that others can verify our results and{/}or adopt the
technique for their own work.

Understandably, there are researchers who 
prefer traditional theoretical approaches.
They may consider computer-assisted proofs as being less than elegant and
even suspicious. We adamantly proclaim that the Sturm procedure is 
mathematically rigorous. If needed, the computation can be written down
line by line in black and white and checked in the traditional theoretical
manner. The computer software is used only so that we can skip these
boring steps, rather than trying to hide  possible 
errors.

Throughout the rest of the paper, we maintain the notation $ I=[0,\pi ] $
and $ J=[-1,1] $.

{\bf Note.}
The author does not intend to publish this paper in a regular journal.
It will be archived in the internet, available to be referenced by
others.

\section{Cosine Polynomials}

\begin{EX} \rm
Let us show that
\begin{equation}  C_1(x) = 5 + 4\,\cos(x) + 3\,\cos(2x) + 4\,\cos(3x) \geq 0  \end{equation}
in $ I $. This assertion does not follow from any known criterion.

For all integers $ n $, $ \cos(nx) $ can be expanded into 
a polynomial of $ \cos(x) $.
For example, $ \cos(3x)=4\,\cos^3(x)-\cos(x)=4y^3-y $, where $ y=\cos(x) $.
Applying this to $ C_1(x) $, we obtain
\begin{equation}  C_1(x) =  P_1(y) = 16y^3 + 6y^2 - 8y + 2 .  \end{equation}
As $ x $ varies from 0 to $ \pi  $, $ y $ varies from $ 1 $ to $ -1 $.
Hence, $ C_1(x)\geq 0 $ in $ I $ if and only if $ P_1(y)\geq 0 $ in $ J=[-1,1] $.
The latter assertion can be easily verified by noting that
$ P_1(y)=2(y+1)(8y^2-5y+1) $, or by
studying the critical points of $ P_1(y) $ in $ J $.
\end{EX}

\begin{EX} \rm
When the given CP is of a higher degree, for instance
\begin{equation}  C_2(x) = 7 + 6\,\cos(x) + 5\,\cos(2x) + 4\,\cos(3x) +3\,\cos(4x) + 5\,\cos(5x) ,  \end{equation}
the resulting AP (algebraic polynomial)
\begin{equation}  P_2(y) = 80y^5 + 24y^4 - 84y^3 - 14y^2 + 19y + 5  \end{equation}
is of higher degree in $ y $. Attempting to verify $ P_2(Y)\geq 0 $ by either
factoring or by examining the critical points may not be feasible in general.
In this example, factoring can simplify the proof a little.
\begin{equation}  P_2(y) = (y+1) X = (y+1)(80y^4-56y^3-28y^2+14y+5) ,  \end{equation}
but the desired assertion $ X\geq 0 $ in $ [-1,1] $ still does not seem obvious.
\end{EX}

We resort to a classical Theorem of Sturm \cite[Chapter 79]{vdw}. Following
van der Waerden, we first differentiate $ X $ to get
\begin{equation}  X_1 = \frac{dX}{dy} = 320y^3-168y^2-56y+14 .  \end{equation}
The classical Euclidean algorithm is then applied to $ X $ and $ X_1 $, to get
a sequence of ``{\em n-remainders}'' $ X_i $ (the so-called Sturm sequence, which
are the negative of the conventional remainders) and quotients $ Q_i: $

\begin{equation}  X_2 = \frac{427}{20} \,y^2 - \frac{161}{20} \, Y - \frac{449}{80} \,, \hspace*{20mmmm} X = Q_1X_1 - X_2  \end{equation}
\begin{equation}  \hspace*{2mmmm}      X_3 = - \frac{267520}{26047} \,y - \frac{40480}{26047} \,, \hspace*{25mmmm} X1 = Q_2X_2 - X_3  \end{equation}
and
\begin{equation}  \hspace*{2.5mmmm}      X_4 = - \frac{360937}{92416} \,, \hspace*{42mmmm} X2 = Q_3X_3 - X_4  \end{equation}
Next we evaluate $ \left\{ X,X_1,. ,X_4\right\}  $ at the two endpoints, respectively, to obtain the
two sequences of numbers:
\begin{equation}  \left\{  99 , -418 , \frac{1903}{80} , \frac{227040}{26047} , \frac{360937}{92416} \right\}   \quad  \mbox{and} \quad       \left\{  15 , 110 , \frac{123}{16} , - \frac{44000}{3721} , \frac{360937}{92416} \right\}  \,.  \end{equation}
In this example, none of the numbers are 0. In the general case, any 0 in the two
sequences are deleted.

In reality, we are only concerned with the signs of these numbers.
There are 2 changes of signs in the first sequence and 2 in the second sequence.
The difference between these two numbers, 0 for the current example,
gives the exact number of real roots 
(each set of multiple roots are counted once only) of $ X $ in $ (-1,1) $.

As a consequence, $ X $ has no roots in $ (-1,1) $ and $ X $ must be of one
sign in the interval. Since $ X(0)=5>0 $, we conclude that
$ X>0 $ in $ (-1,1) $, which in turn implies that 
$ C_2(x)\geq 0 $ in $ I $.

{\bf Remark.} In the above examples, it happens that the two AP, $ P_1(x) $ and 
$ X $, have no roots in $ J $. However, even if the AP has a root in $ J $, as long
as it does not change sign (this occurs when the root is a double root
or a multiple root of even order) in $ J $, we can still conclude that
the corresponding CP is NN in $ I $.

\section{Sine Polynomials}

SP can be dealt with in a similar way. For any integer $ n $, $ \sin(nx) $
can be expanded into a product of $ \sin(x) $ and an AP in $ y=\cos(x) $. For instance,
\begin{equation}  \sin(4x) = \sin(x) \,(8y^3-4y) .  \end{equation}

\begin{EX} \rm
Let us show that
\begin{equation}  S_1(x) =  4\,\sin(x) + 3\,\sin(2x) + 2\,\sin(3x) - 0.8\,\sin(4x) \geq 0  \end{equation}
in $ I $. Again, this is not covered by any known criterion.

After expansion,
\begin{equation}  S_1(x) = \frac{2\sin(x)}{5} \, (16y^3+20y^2+7y+5) .  \end{equation}
Since $ \sin(x)\geq 0 $ in $ I $, $ S_1(x)\geq 0 $ in $ I $ if and only if
\begin{equation}  X = 16y^3+20y^2+7y+5 \geq  0   \end{equation}
in $ J $. Although we do not need Sturm's Theorem to handle this simple
case, we will use it anyway just for illustration. The Sturm sequence for $ X $
is
\begin{equation}  X_1 = 48y^2+40y+7, \quad  X_2=\frac{8}{9} \,y - \frac{145}{36,} \quad  X_3 = - \frac{75123}{64} .  \end{equation}
The two sequences of signs at the two endpoints are the same
\begin{equation}  \left\{ + + - - \right\}  .  \end{equation}
Hence, $ X $ is of one sign in $ J $, implying that $ S_1(x)\geq 0 $ in $ I $.
\end{EX}

\begin{EX} \rm
\begin{equation}  S_2(x) =  8\,\sin(x) + 7\,\sin(2x) + 6\,\sin(3x) + 5\,\sin(4x) + 4\,\sin(5x)\geq 0  \end{equation}
in $ I $.
In a similar way, we can show that $ S_2(x)\geq  $ 0 in $ I $ if and only if 
\begin{equation}  X = 64y^4+40y^3-24y^2-6y+6 \geq 0  \end{equation}
in $ J $. The latter can be rigorously established using Sturm's Theorem.

\section{General Trigonometric Polynomials}
For a general TP that contains both sine and cosine terms, the simple procedure 
described in the previous two sections no longer applies, because after
expansion, we obtain an expression that involves both $ \sin(x) $
and $ \cos(x) $. In some cases, 
a modified procedure may work, as the next example shows. In other cases,
more elaborate arguments are needed. Let us look at a simple example first.
\end{EX}

\begin{EX} \rm
\begin{equation}  T_1(x) = 6 + (6\,\cos(x)+6\,\sin (x))+(-2\,\cos(3x)+6\,\sin(2x))-\cos(4x).  \end{equation}
After expansion, we arrive at the expression
\begin{eqnarray}
   X_0 &=&   (8y^2+12y+4)\,\sin(x)-(8y^4+8y^3-8y^2-12y-5) \nonumber \\
   &=&  (8y^2+12y+4)\,\sqrt{1-y^2} -(8y^4+8y^3-8y^2-12y-5).
\end{eqnarray}
If this has a root in $ J $, then it is also a root of
\begin{equation}  X = (8y^2+12y+4)^2(1-y^2) -(8y^4+8y^3-8y^2-12y-5)^2 = 0,  \end{equation}
obtained by squaring the two parts of $ X_0 $.
The Sturm procedure can be applied to show that $ X $ has no roots in $ J $.
It follows that $ X_0 $ has no root in $ J $
and thus $ T_1(x) $ is of one sign in $ I $. Since $ T_1(0)>0 $, $ T_1(x)\geq 0 $ in $ I $.

A pitfall of this approach is that we only have a one-way
implication (not an ``if and only if'' relation)
between $ X $ having no roots and $ X_0 $ having no roots.
The step of squaring the two parts of $ X_0 $ to obtain $ X $ can introduce
extraneous roots. Thus,
if it happens that $ X $ has roots in $ J $, it is still possible that $ X_0 $
has no roots in $ J $. The next example indicates how such a situation can
be dealt with. The arguments are not particularly elegant; a better
solution will be welcome.
\end{EX}

\begin{EX} \rm
\begin{equation}  T_2(x) = \frac{7}{5} + (\cos(x)+\sin(x))+2\,\sin(2x) + \sin(3x) .  \end{equation}
After expansion, the corresponding $ X_0 $ and $ X $ are given by
\begin{equation}  X_0 =  (4y+4y^2)\,\sqrt{1-y^2} + \left( \frac{7}{5} +y\right)  = f(y)+g(y),  \end{equation}
where $ f(y)=(4y+4y^2)\,\sqrt{1-y^2} $ and $ g(y)=(7/5+y) $ denote the two parts of
$ X_0 $, and
\begin{equation}  X = (4y+4y^2)^2(1-y^2) - \left( \frac{7}{5} +y\right) ^2 .  \end{equation}
Sturm's Theorem reveals that there are two roots of X
in $ J $. Numerical computation locates the roots
in the subintervals $ J_1=[0.345,0.346] $ and $ J_2=[0.948,0.949] $, a fact
that can be rigorously confirmed using Sturm's procedure.
The complement of
$ J_1\cup J_2 $ consists of three subintervals $ [-1,0.345) $, $ (0.346,0.948) $
and $ (0.949,1] $. In each of these
subintervals, $ X $ and, hence, also $ X_0 $ are of one sign, and it is easy to
check that the sign must all be positive. It is also easy to verify that
in $ J_1 $, both $ f(y) $ and g(y) are positive so that $ X_0=f(y)+g(y)\geq 0 $. The same
is true for $ J_2 $. Thus we conclude that $ X_0\geq 0 $ in $ J $, and so $ T_2(x)\geq 0 $
in $ I $.
\end{EX}

\section{The Use of Symbolic Computation Software}

The Sturm procedure described in the previous sections involves a lot
of technical computation. Luckily, 
symbolic computation software can make life easier. We have been using
MAPLE 13 for our investigation.
In this Section, we assume that the reader is familiar with rudimentary
MAPLE syntax and usage, including variables, built-in commands,
functions, procedures, etc. 

\begin{enumerate}
\item Expanding a TP into an expression in $ y $ and $ \sin(x) $ is tedious.
The MAPLE command ``{\tt expand}'' can be used to achieve the goal.

Suppose that the variable {\tt C1} contains the CP as in 
Example 1.
\begin{quote}
\begin{verbatim}
C1 := 5 + 4*cos(x) + 3*cos(2*x) + 4*cos(3*x);
\end{verbatim}
\end{quote}
Then the statement
\begin{quote}
\begin{verbatim}
P1 := subs(cos(x)=y, expand(C1));
\end{verbatim}
\end{quote}
will produce the desired AP and assign it to the variable
{\tt P1}.
On the other hand, to expand an SP, we use a slightly more
complicated statement. Substituting {\tt sin(x)=1} has the same effect
as eliminating the factor $ \sin(x) $.
\begin{quote}
\begin{verbatim}
S1 := 4*sin(x) + 3*sin(2*x) + 2*sin(3*x) - 4*sin(4*x)/5;
X  := subs(sin(x)=1, cos(x)=y, expand(S1));
\end{verbatim}
\end{quote}

Though clumsier, the more complicated statement used for SP works
also for CP.
Since we will be repeating similar computations an enormous number of times,
it is more convenient to define a proc (procedure) 
{\tt pcs} that can be used as a shorthand.
\begin{quote}
\begin{verbatim}
# Define proc pcs(T). T is a variable containing either a sin
#   or cos polynomial. pcs(T) outputs a polynomial in y.
pcs := T -> subs(sin(x)=1, cos(x)=y, expand(T));
P1 := pcs(C1);        # sample usage of pcs
X  := pcs(S1);        # sample usage of pcs
\end{verbatim}
\end{quote}
Though powerful, computers and computer software have their limitations.
In the case of the command {\tt expand}, it knows how to handle 
TP of degree 99 or less, but not higher. An experienced MAPLE user
will no doubt be able to write an extended version of the command to 
jump over the hurdle. We have no need to do that yet.

Another limitation is speed.
As the degree of the CP gets higher, computation takes more time.
The use of computers ought to play a secondary role in theoretical research.
Ideally, in any project, computer assistance should be kept to a minimum,
and used only to take care of a small number of exceptional cases.

\item The MAPLE command ``{\tt sturm(P,y,-1,1)}'' uses the Sturm Theorem to
  count the exact number of real roots of {\tt P} (an AP of the 
  variable {\tt y}) in the interval $ (-1,1] $. MAPLE excludes the left endpoint $ -1 $
  but includes the right endpoint $ 1 $ when counting the roots.

  The official documentation of MAPLE 13, however, describes a different usage.
  It states that one should first use the command ``{\tt sturmseq(P,y)}'' 
  to obtain the Sturm sequence associated with {\tt P} before plugging
  the answer into the command {\tt sturm}. In other words, it recommends
  the clumsier construct {\tt sturm(strumseq(P,y),y,-1,1)}
  for the same goal. The undocumented usage of {\tt sturm}
  described in the previous paragraph was discovered serendipitously
  and has been verified with many test examples.

\item Combining 1 and 2, we have a convenient computer implementation of
  the Sturm procedure. Given any CP or SP {\tt T},
  we simply feed it into the following proc {\tt tsturm}. If the 
  the answer is 0, then the given polynomial is NN in $ I $.
  We do not have an analogous proc for general TP.
\begin{quote}
\begin{verbatim}
# Define proc tsturm(T) that computes the number of roots
#   of a sin/cos poly T in [0,pi]
tsturm := T -> sturm(pcs(T),y,-1,1)
\end{verbatim}
\end{quote}
  Note that the interval $ I $ is hard-wired into the design of the proc
  {\tt tsturm}. To achieve the flexibility of studying more general intervals, 
  the proc can be modified to take the endpoints of the interval as
  input arguments.
\end{enumerate}

\section{Polynomials with Irrational Coefficients}

If all of the coefficients of a TP are rational numbers (and the degree is
less than 99), the simple proc {\tt tsturm} presented in the
previous section is adequate. However, if
some of the coefficients are irrational, the MAPLE command
{\tt sturm} refuses to work. One way to overcome this difficulty is
to write an improved {\tt sturm} command that can handle irrational
coefficients. This can be done but it will involve quite a bit of work.
Instead, we have successfully used an alternate approach 
based on the construction of
a suitable lower bound polynomial with rational coefficients.

\begin{EX} \rm
Consider 
\begin{equation}  S_3(x) =  \left(  \sin(x) + \frac{1}{2} \,\sin(2x) \right)  + \frac{1}{\sqrt{2}} \, \left(  \sin(3x) + \frac{3}{4} \,\sin(4x) \right)  \, .   \Label{s3}  \end{equation}
  It corresponds to the AP
\begin{equation}  P_3(y) = 3\sqrt2 \, y^3 + 2\sqrt2 \, y^2 + \left( 1 -\frac{3}{\sqrt2} \right) \, y +\left( 1-\frac{1}{\sqrt2} \right)  \,.  \end{equation}
  If you feed this to {\tt sturm}, MAPLE complains.

We need to modify the Sturm procedure as follows. 

First we substitute
$ y=z-1 $ to get
\begin{equation}  P_4(z) = P_3(z-1) = 3\sqrt2 \,z^3 - 7\sqrt2 \,z^2 + \left( \frac{7}{\sqrt2} +1 \right) z  .  \Label{p4}  \end{equation}
The reason for doing this is that we prefer to deal with the NN variable 
$ z\in[0,2] $ instead of $ y\in[-1,1] $. In this example,
$ P_4(z) $ is only of degree 3 in $ z $ and, furthermore, 
it has $ z $ as a factor. Thus, we do not need the Sturm procedure
to show that $ P_4(z)\geq 0 $ in $ [0,2] $. Just for the sake of
illustration, let us pretend that we need to use
the Sturm procedure, which cannot be applied directly to $ P_4(z) $.
The next step circumvents the difficulty.

We have devised a MAPLE proc 
``{\tt pfloor(P,m)}'' which takes {\tt P}, an AP in $ z $,
and returns another AP in $ z $, say
{\tt Q}, such that all coefficients of their difference {\tt P-Q} are NN, 
so that {\tt P-Q}$\geq 0$ in $ [0,2] $. The
second argument {\tt m} specifies that each of the coefficients is
less than $ 10^{-m} $. This is used to control the degree of approximation.
A larger {\tt m} gives a more accurate lower bound {\tt Q}.

\begin{quote}
\begin{verbatim}
# Define proc pfloor(p,m) that takes a polynomial p of z and
#   returns a lower bound polynomial with rational coefficients
pfloor := proc(p,m)
   local sc, S, i, pc, pt;
   p := p * 10^m;
   S := 0;
   pc := coeffs(p, z, 'pt');
   for i from 1 to nops([pc]) do
      S := S + floor(evalf(pc[i])-1)*pt[i]; 
   end do;
   S/10^m;
end proc;

Q := pfloor(P,2);      # apply the proc to P
\end{verbatim}
\end{quote}

Suppose now that {\tt P} contains the polynomial given in (\ref{p4}). Then,
after applying {\tt pfloor(P,2)}, we get
\begin{equation}  \mbox{\tt Q} =  \frac{423}{100} \,z^3 - \frac{991}{100} \,z^2 + \frac{593}{100} \,z .  \end{equation}
Now the {\tt sturm} command can be used to verify that $\mbox{\tt Q}\geq 0$
in $ [0,2] $, which in turn
implies that $\mbox{\tt P}\geq 0$ in $ J $ and $ S_3(x)\geq 0 $ in $ I $.
\end{EX}

\begin{EX} \rm
\begin{equation}  S_4(x) =  S_3(x) + \frac{1}{\sqrt{3}} \, \left(  \sin(5x) + \frac{5}{6} \,\sin(6x) \right)  \,    \Label{s4}  \end{equation}
leads to the polynomial
\begin{eqnarray}
   P(z) &=& \frac{80\sqrt3}{9} \,z^5 - \frac{352\sqrt3}{9} \,z^4 + \left(  \frac{176\sqrt3}{3} + 3\sqrt2 \right) z^3  \nonumber \\ && {} - \left(  \frac{308\sqrt3}{9} + 7\sqrt2  \right) z^2 + \left( 1+\frac{55\sqrt3}{9} + \frac{7\sqrt2}{2} \right) z .
\end{eqnarray}
Using the command {\tt Q := pfloor(P,3)} (using $\mbox{\tt m}=2$
yields a less accurate lower bound which has 2 roots in
$[0,2]$) gives
\begin{equation}  Q = \frac{3079}{200} \,z^5 - \frac{8468}{125} \,z^4 + \frac{21171}{200} \,z^3 - \frac{8647}{125} \,z^2 + \frac{16533}{1000} \,z ,  \end{equation}
which has no roots in $ (0,2] $. As a consequence, {\tt P} is NN in $ J $  and $ S_4(x)\geq 0 $
in $ I $.
\end{EX}

\section{Use of the Sturm Procedure in Recent Work}

\subsection{Lower bounds of two TP}

In \cite{akw2}, Alzer and the author prove that for all even integers $ n\geq 2 $
and $ x\in I $,
\begin{equation}  \sum_{k=1}^{n} \frac{\sin(kx)}{k+1} \geq  \frac {(9-\sqrt{137)\sqrt{110-6\sqrt{137}}}}{384} = - 0.0444 ... \, .  \Label{ak9}  \end{equation}
Equality holds if and only if $ n=2 $ and $ x=\cos^{-1}(-(3+\sqrt{137})/16) $.
For odd integers $ n\geq 1 $, the lower bound in (\ref{ak9})
can be improved to 0, a result due
to Brown and Hewitt \cite{bh}. 

(\ref{ak9}) complements a classical result of Rogosinski and Szeg\"o who show that
\begin{equation}  \sum_{k=1}^{n} \frac{\cos(kx)}{k+1} \geq  - \frac{1}{2} .  \Label{rsz}  \end{equation}

The case $ n=2 $ is trivial. For $ n\geq 12 $, we
have a conventional analytic proof. For $ n=4,6,8 $, and $ 10 $, that
proof fails and we need to exploit the Sturm procedure to confirm the inequality.
We use $ n=4 $ as an illustration. The Sturm 
procedure shows that the SP
\begin{equation}  S_5(x) = \sum_{k=1}^{4} \frac{\sin(kx)}{k+1} - \frac{44}{1000} \,\sin(6x)  \Label{s5}  \end{equation}
is NN in $ I $. As a consequence,
\begin{equation}  \sum_{k=1}^{4} \frac{\sin(kx)}{k+1}  \geq  \frac{44}{1000} \,\sin(6x) > - \frac{44}{1000} \,,  \end{equation}
from which (\ref{ak9}), $n=4$, follows. For $ n=6,8 $, and $ 10 $, we use
\begin{equation}  \sum_{k=1}^{n} \,\frac{\sin(kx)}{k+1} - \frac{44}{1000} \,\sin((2n+2)x)  \end{equation}
instead of (\ref{s5}).

It is plausible that if one works hard enough, one may be able to refine
the analytic proof to cover these exceptional cases, but if the 
computer-aided proof works, why waste the extra effort?

In \cite{akw2}, we also prove that for all integers $ n\geq 1 $ and $ x\in I $,
\begin{equation}  \sum_{k=1}^{n} \frac {\sin(kx)+\cos(kx)}{k+1} \geq  - \frac{1}{2} \,,  \end{equation}
another inequality that complements (\ref{rsz}).
Equality holds if and only if $ n=1 $ and $ x=\pi  $.

The odd-degree case is trivial
in view of (\ref{rsz}) and Brown and Hewitt's result. For even $ n $, we
have an analytic proof for $ n\geq 20 $, but have to use the Sturm
procedure to take care of the exceptional cases $ n=2,4,\cdots,18 $.
Since the polynomials under consideration
involve both sine and cosine terms, more elaborate arguments such as
those explained in Section~4 are needed. For an alternate approach, 
see \cite{akw2}.

\subsection{An improved Vietoris sine inequality}

In \cite{kw1}, the following result is established:

An SP $ \,\displaystyle S_6(x)=\sum_{k=1}^{n}a_k\sin(kx) $ is NN in $ I $ if, for $ j=1,2,\cdots $,
\begin{equation}  \frac{(2j-1)\sqrt{j+1}}{2j\sqrt{j}} \,\,a_{2j+1} \, \leq  \, a_{2j}\, \leq \,  \frac{2j-1}{2j} \,\,a_{2j-1} .  \Label{kv}  \end{equation}
This extends the classical Vietoris sine inequality to SP with non-monotone
coefficients. The extreme SP is
\begin{equation}  S_7(x) = \sum_{k=1}^{m} \frac{1}{\sqrt{k}} \left(  \sin((2k-1)x) + \frac{2k-1}{2k} \,\sin(2kx) \right)  \, ,  \Label{s7}  \end{equation}
which includes the examples $ S_3(x) $ of (\ref{s3}) and $ S_4(x) $ of (\ref{s4}) as special
cases. There are also odd-degree sums, obtained by omitting the last 
even-degree term in (\ref{s7}).
The proof of the general result starts with noting that (\ref{s7}) together
with the odd-degree polynomials are NN for $ m\leq 15 $
using the Sturm procedure as described in Section~7. The following is the
MAPLE program used to automate the procedure.

The line numbers are not part of the statements; they are added only
for easy reference.

\begin{verbatim}
 1    pfloor := ... [OMITTED]

 2    psi := k -> sin((2*k-1)*x)+(2*k-1)*sin(2*k*x)/(2*k);

 3    v1 := m -> sum(psi(k)/sqrt(k), k = 1..m);

 4    v2 := proc(n)
 5       if type(n,even) then v1(n/2); 
 6       else v1((n-1)/2) + sin(n*x)/sqrt((n+1)/2);
 7       end if;
 8    end proc;
  
 9    vie := proc(n)
10       local vv, p, q, dp;
11       global acc;
12       vv := v2(n); print(vv); 
13       p := collect(pcs(vv))),y);
13       p := collect(subs(y=z-1,pcs(vv)),z);
14       if type(n,even) then p := collect(normal(p/z),z); end if;
15       print(p);
16       q := pfloor(p, acc);
17       print(q);
18       dp := sort(evalf(p - q));
19       print(dp);
20       if n <> 7 and degree(dp)+1 > nops(dp) then 100; 
21       else sturm(q,z,0,2);
22       end if;
23    end proc;

24    Digits := 30;
25    acc := 9;
26    for i from 3 to 30 do
27       stv := vie(i);
28       if stv > 0 then break; end if;
29    end do;
\end{verbatim}

The definition of the proc {\tt pfloor} has already been described
earlier and is thus omitted.

Line 2 defines an intermediate function {\tt psi(k)} to facilitate the
construction of the SP.

Line 3 defines the function {\tt v1(m)} which is $ S_6(x) $ of degree $ 2m $.

Lines 4 to 8 define the function {\tt v2(n)} which extends {\tt v1}
to include both odd and even-degree SP.

Lines 9 to 23 define a MAPLE proc {\tt vie}. It takes one
argument {\tt n}, which specifies the order of the partial sum. It outputs
either the number of roots of {\tt q} (line 13), or the number 100 (line 12)
if something goes wrong (to be explained below).
Line 12 constructs the SP in question
and assigns it to the variable {\tt vv}, and prints the
result to the screen. Line 13 finds the corresponding AP
and substitutes
{\tt z-1} for {\tt y}. Line 14 gets rid of the factor {\tt z} if {\tt n} is even.
Line 16 finds the polynomial {\tt q} that bounds {\tt p} from below.
Line 18 computes the difference {\tt dp} between the two polynomials {\tt p}
and {\tt q}. It is known that except for one very special case 
($n=7$ and concerning the constant term), all coefficients of
{\tt p} are irrational and hence all coefficients of {\tt dp} must be
strictly positive, unless rounding errors are present; in this 
exceptional case, one or more
terms in {\tt dp} will be missing (because they have zero coefficients). 
Line 20 checks for this possibility, that is,
(except for $ n=7) $ that {\tt nops(dp)}, the number of terms of {\tt dp},
is the same as the degree of the polynomial {\tt p} plus 1.
If this is not true, then the function outputs 100 to alert that something is
wrong. Otherwise, the number of roots are determined in Line 21.

Line 24 sets the floating point computation accuracy to 30 decimal places.

Line 25 sets the approximating accuracy {\tt p - q} to $ 2\cdot10^{-9} $.
Setting it to a lower number causes some {\tt q} to have roots in $ [0,2] $.

Lines 26 to 29 loops from {\tt i = 3} to {\tt 30}. If the loop does not
break prematurely, our claim is confirmed.

There is another place we need the Sturm procedure in the 
proof of the main result for general $ n $. We need
a Lemma which asserts that the two functions
\begin{equation}  \theta _4(x) := \sum_{j=1}^{4}  \left( \frac{1}{\sqrt{k}} - \frac{1}{\sqrt{5}} \right)  \psi _j(x)  \Label{th2}  \end{equation}
\begin{equation}  \theta _{15}(x) := \sum_{j=1}^{15}  \left( \frac{1}{\sqrt{k}} - \frac{1}{4} \right)  \psi _j(x)  \Label{th3}  \end{equation}
are increasing in the intervals $ I_3=[9\pi /64,\pi /2] $, and 
$ I_4=[\pi /2,3\pi /4] $, respectively. 
Here,
\begin{equation}  \psi _j(x) = \sin((2j-1)x) - \frac{2j-1}{2j} \, \sin(2jx) \, .  \Label{F}  \end{equation}
The Lemma is equivalent to the assertion that $ \theta _4'(x) $ and $ \theta _{15}'(x) $
are NN in $ I_3 $ and $ I_4 $, respectively, which can be verified using the Sturm
procedure.

\section{A Related Technique}

\begin{EX} \rm
Suppose we start with a known NN CP or SP, such as
\begin{equation}  2\sin(x) + \sin(2x)+\sin(3x) \geq 0, \quad  x\in I ,  \end{equation}
and vary one of the coefficients to get a more general SP, such as
\begin{equation}  S_8(x) = 2\sin(x) + \sin(2x)+ a \sin(3x) \geq 0, \quad  x\in I .  \end{equation}

Let us ask for what values of the parameter $ a $ will the SP remain NN.
The answer must be a closed interval because the NN property is
preserved under convex combination. Let us determine this interval.
After expanding $ S_8(x) $ and eliminating the factor $ \sin(x) $, we obtain
the AP
\begin{equation}  X = 4ay^2 + 2y +(2-a) .  \end{equation}
Hence, our problem is equivalent to asking for what values of $ a $ will $ X $ have no roots
in $ J $. Since $ X $ is only a quadratic polynomial, it is not difficult to
obtain the exact answer, namely, the closed interval
\begin{equation}  a \in \left[  0 , 1+\sqrt3/2 \right] .  \end{equation}

What is interesting is that the two endpoints of this interval have different
significance with regard to how the NN property is destroyed when $ a $
crosses these boundaries. Plotting the graph of $ X $ ($y\in J$) for
different values of $ a $ will show that as $ a $ decreases from positive to negative
values, $ X $ becomes negative at $ y=-1 $, an endpoint of $J$,
while as $ a $ increases past $ (1+\sqrt3/2) $, $ X $
becomes negative at some interior point of $ J $. When $ a=(1+\sqrt3/2) $, the graph
of $ X $ is tangent to the horizontal axis at some interior point of $ J $.
That interior point is a multiple root of $ X $. The presence of multiple
roots of an AP is indicated by the vanishing of its discriminant.

This observation suggests the following algorithm for finding the
maximal NN interval of the parameter. Substitute the left endpoint of $ J $
into $ X\geq 0 $ and solve for the parameter $ a $ to get a lower bound. 
Do the same for the right
endpoint of $ J $. Then calculate the discriminant of $ X $ as an AP of $ y $,
and find all values of $ a $ that make the discriminant vanish. The endpoints
of the maximal interval are members of the
set of all these computed values of $ a $.
\end{EX}

The discriminant of an AP can be computed using the MAPLE command
{\tt discrim}.

\begin{EX} \rm
Let us add another term to $ S_8(x) $.
\begin{equation}  S_9(x) = 2\sin(x) + \sin(2x)+ a \sin(3x) + \sin(4x) ,  \end{equation}
which leads to the AP
\begin{equation}  X = 8 y^3 + 4ay^2-2y_(2-a) .  \end{equation}
The condition $ X(-1)\geq 0 $ gives $ a>4/3 $, while $ X(1)\geq 0 $ gives no additional 
information. The discriminant of $ X $ is
\begin{equation}  256a^4-512a^3-512a^2+4608a-6656s,  \end{equation}
with roots $ -2 663\ldots $ and $ 1.881648914\ldots $.
It is easy to verify that the
maximal NN interval for $ a $ is $ [4/3,1.881648914\ldots] $.

\end{EX}

\begin{EX} \rm
As another example, consider
\begin{equation}  S_{10}(x) = 36\,\sin(x)+18\,\sin(2x)+28\sin(3x)+21\,\sin(4x)+a\,\sin(5x).  \end{equation}
It is NN in $ I $ when $ a=24 $. It leads to the AP
\begin{equation}  X = 16ay^4+168y^3+(112-12a)y^2-48y+(8+a) .  \end{equation}
$ X(-1)\geq 0 $ implies $ a\geq 0 $.
$ X(1)\geq 0 $ implies $ a\geq -48 $, which is not useful. 
The discriminant
of $ X $, after dividing by 20480, is
\begin{equation}  80\,{a}^{6}-5760\,{a}^{5}+127916\,{a}^{4}-969840\,{a}^{3}+4082209\,{a} ^{2}+30887248\,a-166154688 ,  \end{equation}
which has three positive roots 
$ \displaystyle 4.282\ldots,31.513\ldots,  $
and $ 32.647\ldots $. The Sturm procedure can be used to verify that 
the maximal NN interval is $ [0,31.513\ldots] $. Although $ a=4.282\ldots $
corresponds to a multiple root of $ X $, plotting the graph of $ X $ will show
that that multiple root is actually at $ x\approx -1.45\not\in J $ and the graph of $ X $ in
$ J $ is still well above the horizontal axis.
\end{EX}

\begin{EX} \rm
Here is another variation of the technique. Let us determine what is the maximal
interval of $ a $ so that
\begin{equation}  S_{11}(x) =   \sum_{k=1}^{5} \, (a+5-k)\,\sin(kx) \geq  0 , \quad  \mbox{for } x \in [0,\pi /2] .  \end{equation}
It corresponds to the AP
\begin{equation}  X = 16\, a\, y^4 + (8+8a)\, y^3+(8-8a)\, y^2+(2-2a)\, y+(2+a) .  \end{equation}
$ X(1)\geq 0 $ implies $ a>-4/3 $ and that is the left endpoint of the maximal interval.
The discriminant, after factoring out $ 16384 $, is
\begin{equation}  5\,a^6-160\,a^5+409\,a^4+824\,a^3-449\,a^2+152\,a-25  \end{equation}
and one of its roots $ 28.98537710\ldots $ is the right endpoint of the 
maximal interval. 

If we require $ S_{11}(x)\geq 0 $ in the larger interval $ I $,
then the maximal NN interval is only $ [0,4.1864302648\ldots] $, the right
endpoint of which is another root of the discriminant.
\end{EX}

\end{document}